\documentclass[english]{article}
\usepackage[T1]{fontenc}
\usepackage[latin9]{inputenc}
\usepackage{amsmath}
\usepackage{amssymb}

\makeatletter
\newtheorem{theorem}{Theorem}

\newtheorem{lemma}[theorem]{Lemma}

\newtheorem{proposition}[theorem]{Proposition}

\date{}

\makeatother

\begin{document}

\title{On Decay of Correlations for Exclusion Processes with Asymmetric
Boundary Conditions}

\author{V. A. Malyshev, V. A. Shvets}
\maketitle
\begin{abstract}
We consider a symmetric exclusion process on a discrete interval of
$S$ points with various boundary conditions at the endpoints. We
study the asymptotic decay of correlations as $S\to\infty$. The main
result is asymptotic independence of a stationary distribution whem
the points are far away from each other. We develop a new recurrent
probabilistic approach which is an alternative to Derrida's algebraic
technique. 
\end{abstract}

\section{Introduction}

Exclusion processes are a long-standing popular subject in both mathematics
and theoretical physics. They are a particular case of processes with
local interaction on a lattice. The interest to them is mainly due
to the fact that exclusion processes are the simplest nontrivial model
for collisions in a multiparticle system. Using them, heat conduction
\cite{DELO}, viscosity \cite{LandimOllaVaradhan,MalMan}, quantum
ferromagnet \cite{CaputoMartinelli}, nonequilibrium process \cite{DLS},
etc., models are constructed. Even irrespective of physics, these
processes are natural probabilistic objects. The first systematic
exposition of the relevant theory was given in the well-known monograph
\cite{ligg}.

The simplest exclusion process is a symmetric exclusion process on
the integer lattice $\mathbb{Z}$. This process is well known to have
a continuum of (spatially) uniform invariant measures, which are mixtures
of Bernoulli measures. The same holds for an exclusion process on
a finite segment of the lattice if boundary conditions are empty.
For other boundary conditions, as a rule, it is easy to show that
there are no invariant Bernoulli measures. Using a matrix method (ansatz;
see \cite{DLS,DEHP}), one can obtain an explicit form of correlation
functions for such a process. This would imply that the invariant
measure is asymptotically Bernoullian. It should be noted that this
powerful algebraic method, similar to the famous Bethe ansatz, is
rather cumbersome, is not always mathematically well-grounded, and
has substantial restrictions on its range of applicability. Namely,
we are unaware of any application where jumps of particles can be
longer than 1. It is also important that the probabilistic nature
of the method is absolutely unclear.

In the present paper, we propose another---simple and natural from
the probability theory viewpoint---approach, which can be extended
to jumps of lengths greater than 1, as well as to other boundary conditions.
Being quite different from the Bethe--Derrida methods, our approach
has something in common with them, namely, a certain recursive procedure.
Here we demonstrate the idea of our approach in the simplest situation.
Various generalizations will be considered in further publications.

\section{Problem Statement and the Result}

We consider a simple symmetric exclusion process on an interval $I_{S}=\{0,1,\ldots,S,S+1\}$
of a~1D lattice. The state space for this continuous-time finite
Markov chain \[
(\xi_{t}(0),\xi_{t}(1),\ldots,\xi_{t}(S+1)),\quad\xi_{t}(i)=0,1,\]
 is the set $\{0,1\}^{I_{S}}$. Jumps are defined as follows. For
any time interval $[t,t+dt]$ and any pair $s,s+1$, where $s=0,1,\ldots,S$,
the nodes $s,s+1$ interchange their values with probability $\lambda\, dt$
(independently for different $s$), so that \[
\xi_{t+dt}(s)=\xi_{t}(s+1),\qquad\xi_{t+dt}(s+1)=\xi_{t}(s).\]
 This defines a generator of a Markov process that we call a process
with empty boundary conditions. However, we shall consider boundary
conditions: \[
\xi_{t}(0)\equiv0,\qquad\xi_{t}(S+1)\equiv1.\]
 More precisely, this means that $\xi_{0}(0)=0$ and interchange between
the points $0,1\in I_{S}$ reduces to the following: the pair $\xi_{t}(0)=0$,
$\xi_{t}(1)=1$ becomes $\xi_{t+dt}(0)=0$, $\xi_{t+dt}(1)=0$ with
probability $\lambda\, dt$. The same for the points $S,S+1\in I_{S}$:
$\xi_{0}(S+1)=1$ and the pair $\xi_{t}(S)=0$, $\xi_{t}(S+1)=1$
becomes $\xi_{t+dt}(S)=1$, $\xi_{t+dt}(S+1)=1$ with probability
$\lambda\, dt$.

Denote by $\pi(n_{0},n_{1},\ldots,n_{S+1})=\pi^{(S)}(n_{0},n_{1},\ldots,n_{S+1})$
the stationary measure of this Markov chain. We are interested in
the functions \[
m_{k}^{(S)}(x_{1},x_{2},\ldots,x_{k})=\pi^{(S)}\bigl(n_{x_{1}}=1,\: n_{x_{2}}=1,\:\ldots,\: n_{x_{k}}=1\bigr),\quad x_{1}<\ldots<x_{k}.\]
 For the empty boundary conditions, the process is reducible (because
of the conservation of the number of particles, i.e., the number of
ones), and any stationary measure is a mixture of Bernoulli measures,
i.e., measures that are products of $S+2$ equal measures. As is easy
to show, in the case of our boundary conditions the process is irreducible
(the number of particles is not preserved now), and its invariant
measure is not a product of independent measures. Our goal is to prove
the following correlation decay property for this process.

\begin{theorem} For all $0\leq x\leq S+1$, we have \[
m_{1}^{(S)}(x)=\frac{x}{S+1}.\]

For any $0<x_{1}=x_{1}(S)<x_{2}=x_{2}(S)<S+1$ such that $\frac{x_{k}(s)}{S+1}\to\alpha_{k},k=1,2$
as $S\to\infty$, we have \[
m_{2}^{(S)}(x_{1},x_{2})\to_{S\to\infty}\alpha_{1}\alpha_{2}\]
 \end{theorem}

\section{Proof}

The proof consists of two parts. The first part, reduction to a dual
process (which is easier to study), is rather standard but cumbersome.
The second part, analysis of the dual process, is original.

\subsection{Moment Closeness}

First we obtain equations for moments of the process $\xi_{t}$: \[
m_{k}(t;\{x_{1},\ldots,x_{k}\})=\mathbf{E}\prod_{i=1}^{k}\xi_{t}(x_{i}),\quad x_{1},\ldots,x_{k}\in I_{S},\quad k=1,2,\ldots.\]
 For a process with values 0 and 1, these moments uniquely determine
marginal distributions \[
Pr\Biggl(\,\bigcap_{x\in\mathcal{S}^{\circ}}\{\xi_{t}(x)=\theta_{x}\}\Biggr)\]
 and, of course, vice versa. These equations possess the property
of moment or marginal closeness \cite{IMM}, which makes it possible
to construct a convenient dual process. In what follows, we assume
for simplicity that $\lambda=1$.

We introduce the notion of a cluster. Let a set $Q\subset\mathbb{N}$
be given. A subset $(k_{1}<...<k_{m})=K\subseteq Q$ is called a cluster
if $k_{i+1}-k_{i}=1$ for all $i=1,\ldots,m-1$ and, moreover, $\{k_{1}-1\}\notin Q$
and $\{k_{m}+1\}\notin Q$. In other words, a cluster is a maximal
chain of successive elements of $Q$.

\begin{lemma} Consider a set $X=(x_{1},\ldots,x_{k})$ such that
$0<x_{1}<\ldots<x_{k}<S+1$. Let $X=\bigcup\limits _{j=1}^{p}K_{j}$,
i.e. $X$ is a union of clusters $K_{j}=\bigl(x_{i_{1}}^{(j)},\ldots,x_{i_{q(j)}}^{(j)}\bigr)$.
Then for any $X$ we have the equation \begin{multline*}
\frac{d}{dt}m_{k}(t,\{x_{1},\ldots,x_{k}\})=\sum_{j=1}^{p}m_{k}\bigl(t,\bigl\{ x_{1},\ldots,x_{i_{1}}^{(j)}-1,\ldots,x_{k}\bigr\}\bigr)\\
+\sum_{j=1}^{p}m_{k}\bigl(t,\bigl\{ x_{1},\ldots,x_{i_{q(j)}}^{(j)}+1,\ldots,x_{k}\bigr\}\bigr)-2pm_{k}(t,\{x_{1},\ldots,x_{k}\})\end{multline*}
 with boundary conditions \[
\begin{gathered}m_{k}(t,\{0,x_{2}\ldots,x_{k}\})=0,\\
m_{k}(t,\{x_{1}\ldots,x_{k-1},S+1\})=m_{k-1}(t,\{x_{1}\ldots,x_{k-1}\}).\end{gathered}
\]
 \end{lemma}

Proof. First we derive the equation for a one-point function. The
moment closeness property is based on the linearity of the conditional
probability \[
Pr\bigl(\xi_{t+dt}(x)=1\mid\xi(x-1)=\alpha,\:\xi_{t}(x)=\beta,\:\xi_{t}(x+1)=\gamma\bigr)=\alpha dt+\beta(1-2\, dt)+\gamma\, dt.\]
 By the total probability formula, \[
\begin{aligned}Pr\bigl(\xi_{t+dt}(x)=1\bigr) & =\sum_{\alpha,\beta,\gamma=0}^{1}\bigl(dt(\alpha+\gamma-2\beta)+\beta\bigr)Pr\bigl(\xi_{t}(x-1)=\alpha,\:\xi_{t}(x)=\beta,\:\xi_{t}(x+1)=\gamma\bigr)\\
 & =dt\Bigl(Pr\bigl(\xi_{t}(x-1)=1\bigr)+Pr\bigl(\xi_{t}(x+1)=1\bigr)-2Pr\bigl(\xi_{t}(x)=1\bigr)\Bigr)+Pr\bigl(\xi_{t}(x)=1\bigr),\end{aligned}
\]
 whence we get \[
\frac{d}{dt}m_{1}(t;\{x\})=m_{1}(t;\{x+1\})+m_{1}(t;\{x-1\})-2m_{1}(t;\{x\}).\]

Now we pass to correlation functions of an arbitrary order $k$. Let
$X$ consist of clusters of length~1 only. Since a change during
time $dt$ can only occur at one of the points $x_{i}$, similar computations
for each of them yield \begin{multline*}
\frac{d}{dt}m_{k}(t;\{x_{1},\ldots,x_{k}\})=\sum_{i=1}^{k}\Bigl(m_{k}(t;\{x_{1},\ldots,x_{i}+1,\ldots,x_{k}\})\\[-3pt]
+m_{k}(t;\{x_{1},\ldots,x_{i}+1,\ldots,x_{k}\})\Bigr)-2km_{k}(t;\{x_{1},\ldots,x_{k}\}).\end{multline*}
 Now consider the situation with clusters of larger lengths. As above,
it suffices to consider the case of only one such cluster. Thus, let
$K_{j}=\bigl(x_{i_{1}}^{(j)},\ldots,x_{i_{q(j)}}^{(j)}\bigr)$ be
the $j$th cluster in the partition of~$X$. Changes can only occur
at the endpoints of the cluster. Let us write the general form of
the conditional probability in this case, omitting the index $j$:
\begin{multline*}
Pr\Bigl(\xi_{t+dt}(x_{i_{1}})=1,\:\ldots,\:\xi_{t+dt}(x_{i_{q}})=1\mid\xi_{t}(x_{i_{1}}-1)=\alpha,\\
\xi_{t}(x_{i_{1}})=\beta,\:\xi_{t}(x_{i_{q}})=\gamma,\:\xi_{t}(x_{i_{q}}+1)=\delta\Bigr)=\alpha\gamma\, dt+\beta\gamma(1-2\, dt)+\beta\delta\, dt.\end{multline*}
 Of course, the conditional probability is taken provided that $\xi_{t}(x_{i_{j}})=1,\:\forall j=2,\ldots,q-1$.
Then we obtain \begin{multline*}
Pr\bigl(\xi_{t+dt}(x_{i_{1}})=1,\:\ldots,\:\xi_{t+dt}(x_{i_{q}})=1\bigr)\\
\begin{aligned} & =\sum_{\alpha,\beta,\gamma,\delta=0}^{1}\bigl(\alpha\gamma\, dt+\beta\gamma(1-2\, dt)+\beta\delta\, dt\bigr)Pr\Bigl(\xi_{t}(x_{i_{1}}-1)=\alpha,\:\xi_{t}(x_{i_{1}})=\beta,\\[-7pt]
 & \hskip24ex\xi_{t}(x_{i_{1}}+1)=1,\:\ldots,\:\xi_{t}(x_{i_{q}}-1)=1,\:\xi_{t}(x_{i_{q}})=\gamma,\:\xi_{t}(x_{i_{q}}+1)=\delta\Bigr)\\
 & =dt\Bigl(m_{q}(t;\{x_{i_{1}}-1,\ldots,x_{i_{q}}\})+m_{q}(t;\{x_{i_{1}},\ldots,x_{i_{q}}+1\})\end{aligned}
\\
-2m_{q}(t;\{x_{i_{1}},\ldots,x_{i_{q}}\})\Bigr)+m_{q}(t;\{x_{i_{1}},\ldots,x_{i_{q}}\}),\end{multline*}
 which implies that for an individual cluster we have the equation
\begin{multline*}
\smash[b]{\frac{d}{dt}}m_{q}(t;\{x_{i_{1}},\ldots,x_{i_{q}}\})=m_{q}(t;\{x_{i_{1}}-1,\ldots,x_{i_{q}}\})+m_{q}(t;\{x_{i_{1}},\ldots,x_{i_{q}}+1\})\\
-2m_{q}(t;\{x_{i_{1}},\ldots,x_{i_{q}}\}).\end{multline*}
 Taking the sum over clusters, we obtain the lemma. For convenience,
we write the right-hand side of the equation as follows: \[
\frac{d}{dt}m_{k}(t;\{x_{1},\ldots,x_{k}\})=\sum_{j=1}^{k}\Delta_{x_{j}}^{2}m_{k}(t;\{x_{1},\ldots,x_{k}\}).\]

\subsection{Dual Process}

Denote by $2^{I_{S}}$ the set of (finite) subsets of $I_{S}$, including
the empty set $\emptyset$. Let $(\xi_{t},\: t\ge0)$ be the exclusion
process defined above. A process $(A_{t},\: t\ge0)$ with values in
$2^{I_{S}}$ is said to be \emph{dual\/} to the process $(\xi_{t},\: t\ge0)$
if for all $t\ge0$ we have \[
E\prod_{x\in A_{0}}\xi_{t}(x)=E\prod_{x\in A_{t}}\xi_{0}(x).\]
 In \cite[Section~VIII.1, Theorem~1.1]{ligg} it is proved that a
\emph{symmetric} exclusion process is \emph{self-dual}. Thus, for
instance, analysis of two-point correlation functions can be reduced
to consideration of two particles performing an exclusion walk. To
extend this result to our boundary conditions, consider the above
equations. They lead to the following statement.

\begin{lemma} A process dual to\/ $(\xi_{t},\: t\ge0)$ is an exclusion
process $(A_{t},\: t\ge0)$ for a finite number of particles with
the following modification: 

{[}(a){]} 
\begin{enumerate}
\item \vskip-5pt If one of the particles of $A_{t}$ touches the boundary
$0$, the whole configuration dies (we also say that all particles
reach the absorbing state $0$): \[
\sigma:=\sup\{s>0:\: A_{s}\cap\{0\}\ne\emptyset\}\quad\Longrightarrow\quad A_{\sigma}=\emptyset;\]
 
\item Any particle that reaches the boundary $S+1$ sticks to it. After
that, the other particles walk independently of the stuck particle,
whose coordinate remains to be $S+1$. 
\end{enumerate}
\end{lemma}\vskip-5pt

Note that after one of the particles reaches the right-hand boundary,
other particles either hit~$0$ or reach the point $S+1$. Therefore,
with probability 1, the Markov process thus defined comes in the course
of time to one of the absorbing states: $\emptyset$ or $\{S+1,\ldots,S+1\}$.

As a consequence, we get the following statement.

\begin{proposition} We have \[
m_{k}^{(S)}(A)=\sum_{\ell\in L(A)}Pr(\ell),\qquad|A|=k,\]
 where $L(A)$ is the set of all trajectories of the dual process
that start in $A$ and do not hit the boundary $0$. \end{proposition}

Proof. Let the process $\eta_{t}(x)=\bigl(\eta_{t}^{1}(x_{1}),\ldots,\eta_{t}^{k}(x_{k})\bigr)\in\mathbb{Z}_{+}^{k}$
be an exclusive random walk of $k$ particles in the segment $[0,S+1]$,
where $x=\smash[b]{(x_{1},\ldots,x_{k})\in\mathbb{Z}_{+}^{k}}$ is
the initial disposition of particles. To prove the duality, it suffices
to show that the function $E\prod\limits _{i=1}^{k}\xi_{0}(\eta_{t}^{i}(x_{i}))$
satisfies the same differential equations as the correlation function
\smash{$E\prod\limits _{i=1}^{k}\xi_{t}(x_{i})$}\rule{0pt}{9.5pt}.
By the total probability formula, we have \[
E\prod_{i=1}^{k}\xi_{0}(\eta_{t}^{i}(x_{i}))=\sum_{i_{1}<\ldots<i_{k}}Pr\bigl(\xi_{0}(i_{1})=1,\:\ldots,\:\xi_{0}(i_{k})=1\bigr)Pr\bigl(\eta_{t}^{1}(x_{1})=i_{1},\:\ldots,\:\eta_{t}^{k}(x_{k})=i_{k}\bigr).\]
 Denote \[
P_{t}(i_{1},\ldots,i_{k}):=Pr\bigl(\eta_{t}^{1}(x_{1})=i_{1},\:\ldots,\:\eta_{t}^{k}(x_{k})=i_{n}\bigr).\]
 This function satisfies the same equation as the moments, i.e., \[
\frac{d}{dt}P_{t}(i_{1},\ldots,i_{n})=\sum_{j=1}^{k}\Delta_{i_{j}}^{2}P_{t}(i_{1},\ldots,i_{k}).\]
 Indeed, for a single particle we have \[
\begin{aligned}P_{t+dt}(i) & =\sum_{j=i-1}^{i+1}Pr\bigl(\eta_{t+dt}(x)=i\mid\eta_{t}(x)=j\bigr)Pr\bigl(\eta_{t}(x)=j\bigr)\\
 & =dt\Bigl(Pr\bigl(\eta_{t}(x)=i-1\bigr)+Pr\bigl(\eta_{t}(x)=i+1\bigr)\Bigr)+(1-2\, dt)\P\bigl(\eta_{t}(x)=i\bigr).\end{aligned}
\]
 Hence, \[
\frac{d}{dt}P_{t}(i)=P_{t}(i-1)+P_{t}(i+1)-2P_{t}(i).\]
 Similar arguments also apply to the case of several particles. Thus,
\[
\begin{aligned}\frac{d}{dt}E\prod_{i=1}^{k}\xi_{0}(\eta_{t}^{i}(x_{i})) & =\sum_{i_{1}<\ldots<i_{k}}Pr\bigl(\xi_{0}(i_{1})=1,\:\ldots,\:\xi(i_{k})=1\bigr)\frac{d}{dt}P_{t}(i_{1},\ldots,i_{k})\\
 & =\sum_{i_{1}<\ldots<i_{k}}Pr\bigl(\xi_{0}(i_{1})=1,\:\ldots,\:\xi_{0}(i_{k})=1\bigr)\Biggl\{\sum_{j=1}^{k}\Delta_{i_{j}}^{2}P_{t}(i_{1},\ldots,i_{k})\Biggr\}\\
 & =\sum_{i_{1}<\ldots<i_{k}}\Biggl\{\sum_{j=1}^{k}\Delta_{i_{j}}^{2}Pr\bigl(\xi_{0}(i_{1})=1,\:\ldots,\:\xi_{0}(i_{k})=1)\Biggr\} P_{t}(i_{1},\ldots,i_{k})\\
 & =\sum_{j=1}^{k}\Delta_{x_{j}}^{2}E\prod_{i=1}^{k}\xi_{0}\bigl(\eta_{t}^{i}(x_{i})\bigr).\end{aligned}
\]
 Clearly, boundary conditions are the same as for the original exclusion
process: \[
\begin{gathered}E\xi_{0}\bigl(\eta_{t}^{1}(0)\bigr)\prod_{i=1}^{k-1}\xi_{0}\bigl(\eta_{t}^{i}(x_{i})\bigr)\equiv0,\\
E\prod_{i=1}^{k-1}\xi_{0}\bigl(\eta_{t}^{i}(x_{i})\bigr)\xi_{0}\bigl(\eta_{t}^{i}(S+1)\bigr)=E\prod_{i=1}^{k-1}\xi_{0}\bigl(\eta_{t}^{i}(x_{i})\bigr).\end{gathered}
\]
 Therefore, the process that we have constructed is indeed dual to
the original exclusion process. Summing up and taking into account
the fact that all particles finally get into one of the two absorbing
states, we obtain \[
\begin{aligned}m_{k}(\infty;\{x_{1},\ldots,x_{k}\}) & =E\prod_{i=1}^{k}\xi_{\infty}(x_{i})=E\prod_{i=1}^{k}\xi_{0}\bigl(\eta_{\infty}^{i}(x_{i})\bigr)\\
 & =Pr\bigl(\eta_{\infty}^{1}(x_{1})=S+1,\:\ldots,\:\eta_{\infty}^{k}(x_{k})=S+1\bigr)=P_{\infty}^{(k)},\end{aligned}
\]
 where $P_{\infty}^{(k)}$ is the probability that all particles in
the dual process reach the absorbing state $S+1$.

\subsection{Analysis of the Dual Process}

Consider the exclusion process $\eta_{t}=\eta_{t}^{(\infty)}=(x_{t},y_{t})$
for two particles on the segment $[0,S+1]$. At the initial time instant,
the particles are at points $0<x_{0}(S)<y_{0}(S)<S+1$, respectively,
such that $\frac{x_{0}(S)}{S+1}\to_{S\to\infty}\alpha$ and $\frac{y_{0}(S)}{S+1}\to_{S\to\infty}\beta$.
It is also convenient to assume that the distance between them is
at least 2. Clearly, this does not lose generality. Let us show that
\[
P_{\infty}^{(2)}\to_{S\to\infty}\alpha\beta,\]
 whence the theorem follows.

For an exclusion process $\eta_{t}$, we introduce time intervals
$(T_{k},T_{k}')$ when particles are at distance~$1$ from each other,
assuming that the rest of the time the distances are greater than
$1$. We call $T_{k}$ the~$k$th meeting moment. Using these time
moments, for $k=0,1,\ldots\strut$ we define processes~$\eta^{(k)}(t)$
(which are non-Markovian if $k>0$!). In each of the processes, two
particles with the same initial disposition are walking.

In the process $\eta^{(0)}(t)$, particles walk independently (simple
symmetric random walk) and do not see each other. In the process $\eta^{(k)}(t)$,
the exclusion process works up to time moment $T'_{k}$. After that
the particles become at distance 2 from each other, and from this
point they start walking independently and do not see each other.

Let $P_{k}$ be the probability that in the process $\eta^{(k)}(t)$
both particles come to the point $S+1$ before one of them hits the
origin. Since the number of particle meetings before absorption is
finite with probability $1$, it suffices to prove that for large
$S$ the probabilities $P_{k}$ differ little from $P_{0}$.

First of all, it is easy to prove that asymptotically as $S\to\infty$
we have \[
P_{0}\to\alpha\beta.\]
 Our main method is recurrence relation between the probabilities
$P_{k}$.

First we compare $P_{1}$ and $P_{0}$. Define a random variable $\tau(x,y)$,
the moment of the first meeting of the particles provided that they
started from the points $x$ and $y$. It is clear that for all processes
$\eta^{(k)}(t)$, $k>0$, the distribution of this variable is the
same. Let $x_{\tau}$ and $y_{\tau}=x_{\tau}+1$ be points where the
first and second particles, respectively, are at the first meeting
moment, and let $\P_{0}(x_{\tau}=n)$ be the probability that the
first particle is at point $n$ at this meeting moment. Let $P_{k}(x,y)$
be the probability for two particles in the $k$th problem to reach
$S+1$ if at the initial time moment they are at points $x$ and $y$.
Note also that the event $\bigcup\limits _{n=1}^{S}(x_{\tau}=n)$
comprises also the outcome that the two particles reach $S+1$ earlier
than die, since in this case the distance between them would necessarily
be equal to one, at least at the points $(S,\: S+1)$.

Then, by the total probability formula, \[
P_{1}=\frac{1}{2}\sum_{n=1}^{S-1}P_{0}(x_{\tau}=n)\bigl(P_{0}(n,n+2)+P_{0}(n-1,n+1)\bigr)+\frac{S}{S+1}P_{0}(x_{\tau}=S).\]
 We transform the sum in parentheses, first rewriting each term: \[
\begin{gathered}P_{0}(n,n+2)-P_{0}(n,n+1)=\frac{n(n+2)}{(S+1)^{2}}-\frac{n(n+1)}{(S+1)^{2}}=\frac{n}{(S+1)^{2}},\\
P_{0}(n-1,n+1)-P_{0}(n,n+1)=\frac{(n-1)(n+1)}{(S+1)^{2}}-\frac{n(n+1)}{(S+1)^{2}}=-\frac{n+1}{(S+1)^{2}},\end{gathered}
\]
 and therefore \[
P_{0}(n,n+2)+P_{0}(n-1,n+1)=2P_{0}(n,n+1)-\frac{1}{(S+1)^{2}}.\]
 Thus, we finally obtain an expression for $P_{1}$ via $P_{0}$:
\[
\begin{aligned}P_{1} & =\sum_{n=1}^{S-1}P_{0}(x_{\tau}=n)\biggl(P_{0}(n,n+1)-\frac{1}{2(S+1)^{2}}\biggr)+\frac{S}{S+1}P_{0}(x_{\tau}=S)\\
 & =\sum_{n=1}^{S}P_{0}(x_{\tau}=n)P_{0}(n,n+1)-\frac{1}{2(S+1)^{2}}\sum_{n=1}^{S-1}P_{0}(x_{\tau}=n)\\
 & =P_{0}-\frac{C_{1}}{2(S+1)^{2}}.\end{aligned}
\]
 To get the last equality, we again use the total probability formula.
Note that $C_{1}=\smash[t]{\sum\limits _{n=1}^{S-1}P_{0}(x_{\tau}=n)}$
is a positive constant less than one, since this sum is the probability
that the particles meet at least once before one of them is absorbed
at $0$ or $S+1$.

Now, taking the obtained result into account, we make similar reasoning
for the probability $P_{2}$: \begin{multline*}
\begin{aligned}P_{2} & =\frac{1}{2}\sum_{n=1}^{S-1}P_{0}(x_{\tau}=n)\bigl(P_{1}(n,n+2)+P_{1}(n-1,n+1)\bigr)+\frac{S}{S+1}P_{0}(x_{\tau}=S)\\
 & =\frac{1}{2}\sum_{n=1}^{S-1}P_{0}(x_{\tau}=n)\biggl(P_{0}(n,n+2)+P_{0}(n-1,n+1)\end{aligned}
\\[-3pt]
-\frac{1}{2(S+1)^{2}}\bigl(C_{1}(n,n+2)+C_{1}(n-1,n+1)\bigr)\biggr)+\frac{S}{S+1}P_{0}(x_{\tau}=S),\end{multline*}
 where \[
C_{1}(i,j)=\sum_{n=1}^{S-1}P_{0}\bigl(x_{\tau(i,j)}=n\bigr).\]
 Denote \[
C_{2}=\frac{1}{2}\sum_{n=1}^{S-1}P_{0}(x_{\tau}=n)\left(C_{1}(n,n+2)+C_{1}(n-1,n+1)\right);\]
 then \[
P_{2}=P_{1}-\frac{C_{2}}{2(S+1)^{2}}.\]
 Taking into account that $C_{1}<1$, we roughly estimate $C_{2}$:
\[
C_{2}=\frac{1}{2}\sum_{n=1}^{S-1}P_{0}(x_{\tau}=n)\bigl(C_{1}(n,n+2)+C_{1}(n-1,n+1)\bigr)<\sum_{n=1}^{S-1}\P_{0}(x_{\tau}=n)<1.\]

Analogously, a similar relation can be written for any $k$: \[
P_{k}=P_{k-1}-\frac{C_{k}}{2(S+1)^{2}},\]
 where \[
\begin{gathered}C_{k}(i,j)=\frac{1}{2}\sum_{n=1}^{S-1}P_{0}\bigl(x_{\tau(i,j)}=n\bigr)\bigl(C_{k-1}(n,n+2)+C_{k-1}(n-1,n+1)\bigr),\\
C_{k}:=C_{k}(\alpha,\beta),\quad C_{k}<1.\end{gathered}
\]
 We rewrite this as follows: \[
P_{k}=P_{0}-\frac{1}{2(S+1)^{2}}\sum_{i=1}^{k}C_{i}.\]
 Passing to the limit as $k\to\infty$, we obtain \[
P_{\infty}=P_{0}-\frac{1}{2(S+1)^{2}}\,\smash[t]{\sum_{i=1}^{\infty}C_{i}}.\]
 Thus, we have an expression for the difference of probabilities \[
P_{0}-P_{\infty}=\frac{1}{2(S+1)^{2}}\sum_{i=1}^{\infty}C_{i}.\]
 It remains to estimate this difference and show that it tends to
zero as $S\to\infty$.

By the construction, $C_{k}$ is the probability that two independently
walking particles meet at least $k$ times before one of them is absorbed.
To estimate this quantity, we introduce an auxiliary model.

Consider a process $X_{t}\in[0,S]$, symmetric random walk of a particle
on the segment $[0,S]$ with $X_{0}=1$. Let $\gamma_{k}$ be the
probability for a particle to come to zero $k$ times before reaching
$S$. We claim that $C_{k}$ is not greater than $\gamma_{k}$. Indeed,
$X_{t}+1$ can be interpreted as the distance between two particles
in our original problem at time instant $t$. The condition $X_{0}=1$
means that at the initial time moment we place the particles at the
smallest distance before their meeting, i.e., at distance $2$. If
$X_{t}=S$, one of the particles is absorbed for sure, since the distance
is equal to $S+1$. Therefore, if two independently walking particles
meet $k$ times before one of them is absorbed, then the auxiliary
particle $X_{t}$ hits zero $k$ times before it comes to $S$. Thus,
the probability to hit zero $k$ times in the auxiliary model is (at
least) not less than the probability to meet $k$ times in the original
problem.

Let us calculate $\gamma_{k}$. It is clear that $\smash[t]{\gamma_{1}=\frac{S-1}{S}}$.
If a particle comes to $0$, it makes the next step to the right,
to the point 1, and everything starts over again; hence, \[
\gamma_{k}=\gamma_{k-1}\frac{S-1}{S}.\]
 Finally, we obtain \[
\gamma_{k}=\left(\frac{S-1}{S}\right)^{k}.\]

Now we can again estimate the difference between $P_{0}$ and $P_{\infty}$:
\[
\begin{aligned}P_{0}-P_{\infty} & =\frac{1}{2(S+1)^{2}}\sum_{i=1}^{\infty}C_{i}\le\frac{1}{2(S+1)^{2}}\sum_{i=1}^{\infty}\gamma_{i}\\
 & =\frac{1}{2(S+1)^{2}}\sum_{i=1}^{\infty}\left(\frac{S-1}{S}\right)^{i}=\frac{S-1}{2(S+1)^{2}}=\frac{1}{2(S+1)}-\frac{1}{(S+1)^{2}}.\end{aligned}
\]
 Therefore, \smash{[}b{]}{$|P_{\infty}-P_{0}|=O\Bigl(\frac{1}{S+1}\Bigr)$}
as $S\to\infty$; since \smash{[}b{]}{$P_{k}\to_{k\to\infty}P_{\infty}$}
for any $S$, this completes the proof.


\begin{thebibliography}{8}
\bibitem{DELO} Derrida, B., Enaud, C., Landim, C., and Olla, S.,
Fluctuations in the Weakly Asymmetric Exclusion Process with Open
Boundary Conditions, \emph{J.~Statist.\ Phys.}, 2005, vol.~118,
no.~5--6, pp.~795--811.

\bibitem{LandimOllaVaradhan} Landim, C., Olla, S., and Varadhan,
S.R.S., On Viscosity and Fluctuation-Dissipation in Exclusion Processes,
\emph{J.~Statist.\ Phys.}, 2004, vol.~115, no.~1--2, pp.~323--363.

\bibitem{MalMan} Malyshev, V. and Manita, A., Stochastic Micro-model
of the Couette Flow, to appear in \emph{Teor.\ Verojatnost.\ i Primenen}.

\bibitem[4]{CaputoMartinelli} Caputo, P. and Martinelli, F., Relaxation
Time of Anisotropic Simple Exclusion Processes and Quantum Heisenberg
Models, \emph{Ann.\ Appl.\ Probab.}, 2003, vol.~13, no.~2, pp.~691--721.

\bibitem[5]{DLS} Derrida, B., Lebowitz, J.L., and Speer, E.R., Entropy
of Open Lattice Systems, \emph{J.~Statist.\ Phys.}, 2007, vol.~126,
no.~4--5, pp.~1083--1108.

\bibitem[6]{ligg} Liggett, T.M., \emph{Interacting Particle Systems},
New York: Springer, 1985. Translated under the title \emph{Markovskie
protsessy s lokal'nym vzaimodeistviem}, Moscow: Mir, 1989.

\bibitem[7]{DEHP} Derrida, B., Evans, M.R., Hakim, V., and Pasquier,
V., Exact Solution of a 1D Asymmetric Exclusion Model Using a Matrix
Formulation, \emph{J.~Phys.~A}, 1993, vol.~26, no.~7, pp.~1493--1517.

\bibitem[8]{IMM} Ignatyuk, I.A., Malyshev, V.A., and Molchanov, S.A.,
Moment-Closed Processes with Local Interaction, \emph{Selecta Math.\ Soviet.},
1989, vol.~8, no.~4, pp.~351--384.
\end{thebibliography}
\end{document}